\documentclass[a4paper,11pt]{amsart}

\usepackage{graphicx}
\usepackage{mathptmx}
\usepackage{amsmath}
\usepackage{amssymb}
\usepackage{enumitem}
\usepackage{xcolor}
\usepackage{pgfplots}

\newmuskip\pFqmuskip

\newcommand*\pFq[6][8]{%
  \begingroup 
  \pFqmuskip=#1mu\relax
  \mathcode`=\string"8000
  \begingroup\lccode`\~=`\,
  \lowercase{\endgroup\let~}\pFqcomma
  F^{#2}_{#3}{\left(\genfrac..{0pt}{}{#4}{#5}\bigg|#6\right)}%
  \endgroup
}
\newcommand{\pFqcomma}{\mskip\pFqmuskip}

\newtheorem{theorem}{Theorem}[section]
\newtheorem{lemma}[theorem]{Lemma}

\newtheorem{remark}[theorem]{Remark}

\begin{document}

\title[Degenerate Euler-Seidel Matrix Method and Their Applications]{Degenerate Euler-Seidel Matrix Method and Their Applications}

\author{Taekyun  Kim}
\address{Department of Mathematics, Kwangwoon University, Seoul 139-701, Republic of Korea}
\email{tkkim@kw.ac.kr}
\author{Dae San  Kim}
\address{Department of Mathematics, Sogang University, Seoul 121-742, Republic of Korea}
\email{dskim@sogang.ac.kr}


\maketitle

\begin{abstract}
This paper introduces a degenerate version of the Euler-Seidel matrix method by incorporating a parameter $\lambda$ into the classical recurrence relation. The standard Euler-Seidel method relates the generating functions of an initial sequence  and its final sequence via Seidel's formula, $\overline{A}(t)=e^{t}A(t)$. Our generalized method establishes transformation formulas using $\lambda$-generalized binomial identities and yields a degenerate Seidel's formula for the exponential generating functions: $\overline{S_{\lambda}}(t)=e_{\lambda}^{1-\lambda}(t)S_{\lambda}(t)$. The results are applied to study and derive new combinatorial identities for sequences like the degenerate Bell and Fubini numbers and polynomials.
\end{abstract}

\section{Introduction}
For any nonzero $\lambda\in\mathbb{R}$, the degenerate exponentials are defined by Kim-Kim as
\begin{equation}
e_{\lambda}^{x}(t)=\sum_{k=0}^{\infty} (x)_{k,\lambda}\frac{t^{k}}{k!},\quad e_{\lambda}(t)=e_{\lambda}^{1}(t), \quad (\mathrm{see}\ [12-17]), \label{1}
\end{equation}
where the generalized falling factorial sequence is given by
\begin{equation}
(x)_{0,\lambda}=1,\quad (x)_{n,\lambda}=x(x-\lambda)\cdots\big(x-(n-1)\lambda\big),\ (n\ge 1). \label{2}
\end{equation}
Also, the generalized rising factorial sequence is defined by
 \begin{equation}
 \langle x\rangle_{0,\lambda}=1,\quad \langle x\rangle_{n,\lambda}=x(x+\lambda)(x+2\lambda)\cdots\big(x+(n-1)\lambda\big),\ (n\ge 1).\label{3}
 \end{equation}
The degenerate Stirling numbers of the second kind are given by
 \begin{equation}
 (x)_{n,\lambda}=\sum_{k=0}^{n}{n \brace k}_{\lambda}(x)_{k},\quad (n\ge 0),\quad (\mathrm{see}\ [13-15]),\label{4}
 \end{equation}
 where $(x)_{0}=1,\ (x)_{n}=x(x-1)\cdots(x-n+1),\ (n\ge 1)$. \\
Then, by \eqref{4}, we get
 \begin{equation}
 \frac{1}{k!}\big(e_{\lambda}(t)-1\big)^{k}=\sum_{n=k}^{\infty}{n \brace k}_{\lambda}\frac{t^{n}}{n!},\quad (k\ge 0),\quad (\mathrm{see}\ [13-15]). \label{5}
 \end{equation}
Note from \eqref{4} or \eqref{5} that
 \begin{displaymath}
 \lim_{\lambda\rightarrow 0}{n \brace k}_{\lambda}={n \brace k},
 \end{displaymath}
where ${n \brace k}$ are the classical Stirling numbers of the second kind defined by
\begin{displaymath}
x^{n}=\sum_{k=0}^{n}{n \brace k}(x)_{k},\quad (n\ge 0),\quad (\mathrm{see}\ [6,24-26]).
\end{displaymath}
The degenerate Bell polynomials are given by
\begin{equation}
\phi_{n,\lambda}(x)=\sum_{k=0}^{n}{n \brace k}_{\lambda}x^{k},\quad (n\ge 0),\quad (\mathrm{see}\ [14,16,17]).\label{6}
\end{equation}
When $x=1, \phi_{n,\lambda}=\phi_{n,\lambda}(1)$ are called the degenerate Bell numbers. \par
From \eqref{1} and \eqref{6}, we have
\begin{equation}
e^{x(e_{\lambda}(t)-1)}=\sum_{n=0}^{\infty}\phi_{n,\lambda}(x)\frac{t^{n}}{n!},\quad (\mathrm{see}\ [14,16,17]). \label{7}
\end{equation}
Note that
\begin{displaymath}
\lim_{\lambda\rightarrow 0}\phi_{n,\lambda}(x)=\phi_{n}(x),
 \end{displaymath}
where $\phi_{n}(x)=\sum_{k=0}^{n}{n \brace k}x^{k}$ are the ordinary Bell polynomials (see [3,4,5,10]). \\
The Fubini polynomials are defined by
\begin{equation}
F_{n}(x)=\sum_{k=0}^{n}{n \brace k}k!x^{k},\quad (n\ge 0),\quad (\mathrm{see}\ [11,20,23,27-29]).\label{8}
\end{equation}
In view of \eqref{6} and \eqref{8}, the degenerate Fubini polynomials are given by
\begin{equation}
F_{n,\lambda}(x)=\sum_{k=0}^{n}{n \brace k}_{\lambda}k!x^{k},\quad (n\ge 0),\quad (\mathrm{see}\ [20,29]). \label{9}
\end{equation}
Thus, by \eqref{1} and \eqref{9}, we get
\begin{equation}
\frac{1}{1-x(e_{\lambda}(t)-1)}=\sum_{n=0}^{\infty}F_{n,\lambda}(x)\frac{t^{n}}{n!},\quad (\mathrm{see}\ [20,29]). \label{10}	
\end{equation}
When $x=1$, $F_{n,\lambda}=F_{n,\lambda}(1)$ are called the degenerate Fubini numbers (or the degenerate ordered Bell numbers). As general references for this paper, the reader may refer to [1,6,24-26].\par
\vspace{0.1in}
For a given sequence $(a_{n})$, the Euler-Seidel matrix corresponding to this sequence is determined recursively by
\begin{equation}
a_{n,0}=a_{n},\ (n\ge 0),\quad a_{n,k}=a_{n,k-1}+a_{n+1,k-1},\quad (n\ge 0,\, k\ge 1). \label{11}
\end{equation}
The Euler-Seidel matrix $(a_{k,n})_{n,k \ge 0}$ corresponding to the sequence $(a_{n})$ is given by
 \begin{equation}
 A=\left(\begin{matrix}
 a_{0,0} & a_{1,0} & a_{2,0} & a_{3,0} & a_{4,0} & \cdots\\
 a_{0,1} & a_{1,1} & a_{2,1} & a_{3,1} & a_{4,1} & \cdots\\
 a_{0,2} & a_{1,2} & a_{2,2} & a_{3,2} & a_{4,2} & \cdots\\
 \vdots & \vdots & \vdots & \vdots & \vdots & \vdots
 \end{matrix}\right).\label{12}
 \end{equation}
The sequence $(a_{n,0})_{n \ge 0}$ is the initial sequence, and $(a_{0,n})_{n\ge 0}$ is the final sequence. \\
Then, by \eqref{11}, we easily get the binomial identities:
\begin{equation}
a_{0,n}=\sum_{k=0}^{n}\binom{n}{k}a_{k,0},\quad a_{n,0}=\sum_{k=0}^{n}\binom{n}{k}(-1)^{n-k}a_{0,k}, \label{13}
\end{equation}
where $n$ is nonnegative integer (see [2,7,8,9,21,22]). \\
Let $A_{1}(t)=\sum_{n=0}^{\infty}a_{n,0}t^{n}$ be the generating function of the initial sequence $(a_{n,0})_{n\ge 0}$. Then the generating function of the final sequence $(a_{0,n})_{n\ge 0}$, due to Euler, is given by
\begin{equation}
\overline{A_{1}}(t)=\sum_{n=0}^{\infty}a_{0,n}t^{n}=\frac{1}{1-t}A_{1}\bigg(\frac{t}{1-t}\bigg),\quad (\mathrm{see}\ [9]).\label{14}
\end{equation}
Let $A_{2}(t)=\sum_{n=0}^{\infty}a_{n,0}\frac{t^{n}}{n!}$ be the exponential generating function of the initial sequence $(a_{n,0})_{n\ge 0}$. Then the exponential generating function of the final sequence $(a_{0,n})_{n\ge 0}$, due to Seidel, is given by (see \eqref{14})
\begin{equation}
\overline{A_{2}}(t)=\sum_{n=0}^{\infty}a_{0,n}\frac{t^{n}}{n!}=e^{t}A_{2}(t),\quad (\mathrm{see}\ [9]).\label{15}
\end{equation} \par
The aim of this paper is to extend these results by studying a degenerate version of the Euler-Seidel method. This generalization introduces a parameter $\lambda$ into the recurrence relation.
For a given sequence $(a_{n,\lambda})$, the {\it{degenerate Euler-Seidel matrix}} is recursively defined by (see \eqref{11}):
\begin{equation*}
\begin{aligned}
&a_{n,0}(\lambda)=a_{n,\lambda},\quad (n\ge 0), \\
&a_{n,k}(\lambda)=\big(1-(k-n)\lambda\big)a_{n,k-1}(\lambda)+a_{n+1,k-1}(\lambda),\quad (n\ge 0,\, k\ge 1).
\end{aligned}
\end{equation*}
The degenerate Euler-Seidel matrix $(a_{k,n}(\lambda))_{n,k \ge 0}$ corresponding to the sequence $(a_{n,\lambda})$ is given by (see \eqref{12})
\begin{equation*}
A(\lambda)=\left(\begin{matrix}
a_{0,0}(\lambda) & a_{1,0}(\lambda) & a_{2,0}(\lambda) & a_{3,0}(\lambda) & \cdots \\
a_{0,1}(\lambda) & a_{1,1}(\lambda) & a_{2,1}(\lambda) & a_{3,1}(\lambda) & \cdots \\
a_{0,2}(\lambda) & a_{1,2}(\lambda) & a_{2,2}(\lambda) & a_{3,2}(\lambda) & \cdots \\
\vdots & \vdots & \vdots & \vdots & \vdots
\end{matrix}\right).
\end{equation*}
The sequence $\big(a_{n,0}(\lambda)\big)_{n \ge 0}$ is the {\it{initial degenerate sequence}}, and $\big(a_{0,n}(\lambda)\big)_{n\ge 0}$ is the {\it{final degenerate sequence}}. The following $\lambda$-generalized binomial identities are established (see Theorems 2.2, 2.3) using the generalized falling and rising factorials, $(1-\lambda)_{n-k,\lambda}$ and $\langle1-\lambda \rangle_{n-k,\lambda}$, (see \eqref{2}, \eqref{3}) for the degenerate case:
\begin{align*}
&a_{0,n}(\lambda)=\sum_{k=0}^{n}\binom{n}{k}(1-\lambda)_{n-k,\lambda}a_{k,0}(\lambda), \\
&a_{n,0}(\lambda)= \sum_{k=0}^{n}\binom{n}{k}(-1)^{n-k} \langle1-\lambda \rangle_{n-k,\lambda}a_{0,k}(\lambda).
\end{align*}
These correspond directly to the classical binomial identities in \eqref{13}, and they lead to a degenerate version of Seidel's formula (see \eqref{15}) for the exponential generating functions
\begin{equation*}
\overline{S_{\lambda}}(t)=\sum_{n=0}^{\infty}a_{0,n}(\lambda)\frac{t^{n}}{n!}=e_{\lambda}^{1-\lambda}(t)S_{\lambda}(t),
\end{equation*}
where $S_{\lambda}(t)=\sum_{n=0}^{\infty}a_{n,0}(\lambda)\frac{t^{n}}{n!}$ is the exponential generating function of the initial degenerate sequence (see Theorem 2.4).
The results derived from the degenerate Euler-Seidel method are applied to well-known sequences, including the degenerate Bell numbers and polynomials, and degenerate Fubini numbers and polynomials, yielding various combinatorial identities (see Theorems 2.5-2.14).\par
The study of degenerate versions of certain special polynomials and numbers has recently seen renewed interest among mathematicians. They have been explored by employing various tools, including generating functions, combinatorial methods, umbral calculus, operator theory, $p$-adic analysis, probability theory, special functions, differential equations and analytic number theory (see [12-19,20,29] and the references therein).
\vspace{0.1in} \par
The following facts will be repeatedly used throughout this paper.
\begin{lemma}
For any nonnegative integer $n$, the following hold true.
\begin{flalign*}
&(a)\,\, (-x)_{n,\lambda}=(-1)^{n}\langle x \rangle_{n,\lambda},\\
&(b)\,\, \langle -x \rangle_{n,\lambda}=(-1)^{n}(x)_{n,\lambda}, \\
&(c)\,\,(x+y)_{n,\lambda}=\sum_{k=0}^{n}\binom{n}{k}(x)_{k,\lambda}(y)_{n-k,\lambda},\\
&(d)\,\, \sum_{k=j}^{n}\binom{n}{k}(x)_{n-k,\lambda}\binom{k}{j}(y)_{k-j,\lambda}=\binom{n}{j}(x+y)_{n-j,\lambda},\\
&(e)\,\,\binom{x+1}{n}=\binom{x}{n}+\binom{x}{n-1}. &&
\end{flalign*}
\begin{proof}
(a) and (b) follow from definitions (see \eqref{2}, \eqref{3}). (c) is immediate from $e_{\lambda}^{x+y}(t)=e_{\lambda}^{x}(t)e_{\lambda}^{y}(t)$ (see \eqref{1}). (d) follows using (c) and $\binom{n}{k}\binom{k}{j}=\binom{n}{j}\binom{n-j}{k-j}$. (e) is well known and easy to show.
\end{proof}
\end{lemma}

\section{Degenerate Euler-Seidel Matrix Method and Their Applications}
Given a sequence $(a_{n,\lambda})$, we consider the degenerate Euler-Seidel matrix corresponding to this sequence which is determined recursively by
\begin{equation}
\begin{aligned}
&a_{n,0}(\lambda)=a_{n,\lambda},\quad (n\ge 0), \\
&a_{n,k}(\lambda)=\big(1-(k-n)\lambda\big)a_{n,k-1}(\lambda)+a_{n+1,k-1}(\lambda),\quad (n\ge 0,\, k\ge 1).
\end{aligned}\label{16}
\end{equation}

\begin{lemma}
For any integers $l$ with $0 \le l \le n$, the following holds true.
\begin{equation}
a_{1,n}(\lambda)=\sum_{k=0}^{l}\binom{l}{k}\big(1-(n-l)\lambda\big)_{l-k,\lambda}a_{k+1,n-l}(\lambda). \label{17}
\end{equation}
\begin{proof}
The identity holds true for $l=0, 1$. Assume that it holds for $l$ with $1 \le l \le n-1$. Then, by assumption and \eqref{16}, we have
\begin{align*}
a_{1,n}(\lambda)&=\sum_{k=0}^{l}\binom{l}{k}\big(1-(n-l)\lambda\big)_{l-k,\lambda}a_{k+1,n-l}(\lambda) \\
&=\sum_{k=0}^{l}\binom{l}{k}\big(1-(n-l)\lambda\big)_{l-k,\lambda} \\
&\quad\times \Big(\big(1-(n-l-k-1)\lambda \big)a_{k+1,\, n-l-1}(\lambda)+a_{k+2,\,n-l-1}(\lambda) \Big)\\
&=\sum_{k=0}^{l}\binom{l}{k}\big(1-(n-l)\lambda\big)_{l-k,\lambda}\big(1-(n-l-1)\lambda \big)a_{k+1,n-l-1}(\lambda)
\end{align*}
\begin{align*}
&\quad +\sum_{k=0}^{l}\binom{l}{k} k\lambda \big(1-(n-l)\lambda\big)_{l-k,\lambda} a_{k+1,\,n-l-1}(\lambda) \\
& \quad +\sum_{k=0}^{l}\binom{l}{k}\big(1-(n-l)\lambda\big)_{l-k,\lambda}a_{k+2,\,n-l-1}(\lambda) \\
&=\sum_{k=0}^{l+1}\binom{l}{k}\big(1-(n-l-1)\lambda\big)_{l-k+1,\,\lambda}a_{k+1,\, n-l-1}(\lambda) \\
&\quad +\sum_{k=0}^{l+1}\bigg\{\binom{l}{k}k \lambda+\binom{l}{k-1}\big(1-(n-k)\lambda\big)\bigg\}\\
&\quad \times \big(1-(n-l)\lambda \big)_{l-k,\lambda} a_{k+1,\, n-l-1}(\lambda) \\
&=\sum_{k=0}^{l+1}\binom{l}{k}\big(1-(n-l-1)\lambda\big)_{l-k+1,\,\lambda}a_{k+1,\, n-l-1}(\lambda) \\
&\quad +\sum_{k=0}^{l+1}\binom{l}{k-1}\big(1-(n-l-1)\lambda\big)_{l-k+1,\lambda}a_{k+1,\,n-l-1}(\lambda) \\
&=\sum_{k=0}^{l+1}\binom{l+1}{k}\big(1-(n-l-1)\lambda \big)_{l+1-k,\lambda}a_{k+1,\,n-l-1}(\lambda),
\end{align*}
which shows \eqref{17} holds for $l+1$.
\end{proof}
\end{lemma}
By using Lemma 2.1, we prove the following theorem.
\begin{theorem}
For $n\ge 0$, we have
\begin{equation}
a_{0,n}(\lambda)=\sum_{k=0}^{n}\binom{n}{k}(1-\lambda)_{n-k,\lambda}a_{k,0}(\lambda). \label{18}
\end{equation}
\begin{proof}
We show this by induction on $n$.
The identity holds trivially for $n=0$. Assume that \eqref{18} holds for $n \ge 0$. Then, by assumption and using \eqref{16} and \eqref{17} with $l=n$, we have
\begin{align*}
&a_{0,n+1}(\lambda)=\big(1-(n+1)\lambda \big)a_{0,n}(\lambda)+a_{1,n}(\lambda) \\
&=\big(1-(n+1)\lambda \big)\sum_{k=0}^{n}\binom{n}{k}(1-\lambda)_{n-k,\lambda}a_{k,0}(\lambda)
+\sum_{k=0}^{n}\binom{n}{k}(1-\lambda)_{n-k-1,\lambda}a_{k+1,0}(\lambda) \\
&=\sum_{k=0}^{n}\binom{n}{k}\big(1-(n-k+1)\lambda \big)(1-\lambda)_{n-k,\lambda}a_{k,0}(\lambda) \\
&\quad +\sum_{k=0}^{n}\binom{n}{k}(-k \lambda)(1-\lambda)_{n-k,\lambda}a_{k,0}(\lambda)
+\sum_{k=0}^{n}\binom{n}{k}(1-\lambda)_{n-k-1,\lambda}a_{k+1,0}(\lambda) \\
&=\sum_{k=0}^{n+1}\binom{n}{k}(1-\lambda)_{n-k+1,\lambda}a_{k,0}(\lambda)
+\sum_{k=0}^{n+1}\bigg\{-\lambda\binom{n}{k}k+\binom{n}{k-1}\bigg\}(1-\lambda)_{n-k,\lambda}a_{k,0}(\lambda) \\
&=\sum_{k=0}^{n+1}\binom{n}{k}(1-\lambda)_{n-k+1,\lambda}a_{k,0}(\lambda)
+\sum_{k=0}^{n+1}\binom{n}{k-1}(1-\lambda)_{n-k+1,\lambda}a_{k,0}(\lambda)\\
&=\sum_{k=0}^{n+1}\binom{n+1}{k}(1-\lambda)_{n-k+1,\lambda}a_{k,0}(\lambda),
\end{align*}
which completes the proof.
\end{proof}
\end{theorem}
For $n\ge 0$, by Lemma 1.1, we get
\begin{align}
&\sum_{k=0}^{n}\binom{n}{k}(-1)^{n-k}\langle 1-\lambda\rangle_{n-k,\lambda}a_{0,k}(\lambda)\label{19}\\
&=\sum_{k=0}^{n}\binom{n}{k}(-1)^{n-k}\langle 1-\lambda\rangle_{n-k,\lambda}\sum_{j=0}^{k}\binom{k}{j}(1-\lambda)_{k-j,\lambda}a_{j,0}(\lambda)\nonumber\\
&=\sum_{j=0}^{n}a_{j,0}(\lambda)\sum_{k=j}^{n}\binom{n}{k}	(\lambda-1)_{n-k,\lambda}\binom{k}{j}(1-\lambda)_{k-j,\lambda} \nonumber\\
&=\sum_{j=0}^{n}a_{j,0}(\lambda)\binom{n}{j}(0)_{n-j,\lambda}=a_{n,0}(\lambda). \nonumber
\end{align}
Therefore, by \eqref{19}, we obtain the inverse relation of \eqref{18}.
\begin{theorem}
For $n\ge 0$, we have
\begin{equation}
a_{n,0}(\lambda)= \sum_{k=0}^{n}\binom{n}{k}(-1)^{n-k}\langle 1-\lambda\rangle_{n-k,\lambda}a_{0,k}(\lambda).\label{20}
\end{equation}
\end{theorem}

Let
\begin{displaymath}
S_{\lambda}(t)=\sum_{n=0}^{\infty}a_{n,0}(\lambda)\frac{t^{n}}{n!}
\end{displaymath}
be the generating function of the initial degenerate sequence $\big(a_{n,0}(\lambda)\big)_{n\ge 0}$.\par
Then, by Theorem 2.2, we have
\begin{align}
e_{\lambda}^{1-\lambda}(t)S_{\lambda}(t)&=\sum_{l=0}^{\infty}(1-\lambda)_{l,\lambda}\frac{t^{l}}{l!}\sum_{k=0}^{\infty}a_{k,0}(\lambda)\frac{t^{k}}{k!} \label{21}\\
&=\sum_{n=0}^{\infty}\sum_{k=0}^{n}\binom{n}{k}(1-\lambda)_{n-k,\lambda}a_{k,0}(\lambda)\frac{t^{k}}{k!}\nonumber\\
&=\sum_{n=0}^{\infty}a_{0,n}(\lambda)\frac{t^{n}}{n!}.\nonumber	
\end{align}
Therefore, by \eqref{21}, we obtain the {\it{degenerate Seidel's formula}}.
\begin{theorem}
Let $S_{\lambda}(t)=\sum_{n=0}^{\infty}a_{n,0}(\lambda)\frac{t^{n}}{n!}$ be the generating function of the initial degenerate sequence $(a_{n,0}(\lambda))_{n\ge 0}$. Then the generating function of the final degenerate sequence $(a_{0,n}(\lambda))_{n\ge 0}$ is given by
\begin{equation}
\overline{S_{\lambda}}(t)=\sum_{n=0}^{\infty}a_{0,n}(\lambda)\frac{t^{n}}{n!}=e_{\lambda}^{1-\lambda}(t)S_{\lambda}(t). \label{22}
\end{equation}
\end{theorem}
\begin{remark}
Assume that $a_{n,\lambda} \rightarrow a_{n}$, as $\lambda \rightarrow 0$, for all $n \ge 0$. Then, from \eqref{11} and \eqref{16}, we see that $a_{n,k}(\lambda) \rightarrow a_{n,k}$, as $\lambda \rightarrow 0$, for all $n,\,k \ge 0$. This is because the recurrence relations in \eqref{16} converge to those ones in \eqref{11}, as $\lambda \rightarrow 0$. In particular, the degenerate Seidel's formula in \eqref{22} becomes the original Seidel's formula in \eqref{15}. Thus our approach not only preserves the structure of the classical Euler-Seidel method as the parameter $\lambda \to 0$ but also provides a powerful framework for studying degenerate versions of combinatorial sequences.

\end{remark}
For example, let $(a_{n,0}(\lambda))_{n\ge 0}=(\phi_{n,\lambda})_{n\ge 0}$, (see \eqref{6}). Then degenerate Euler-Seidel matrix corresponding to this sequence is given by
\begin{equation*}
\left(\begin{matrix}
1 & 1 & 2-\lambda & \cdots  \\
2-3\lambda & 3-\lambda & 7-5\lambda+\lambda^{2} & \cdots \\
5-6\lambda+2\lambda^{2} & 10-9\lambda+2\lambda^{2} & 27-37\lambda-19\lambda^{2}+16\lambda^{3}  & \cdots \\
\vdots & \vdots & \vdots  & \vdots
\end{matrix}\right).
\end{equation*}
From \eqref{7}, we note that
\begin{align}
S_{\lambda}(t)&=\sum_{n=0}^{\infty}a_{n,0}(\lambda)\frac{t^{n}}{n!}=\sum_{n=0}^{\infty}\phi_{n,\lambda}\frac{t^{n}}{n!}=e^{e_{\lambda}(t)-1}.\label{23}
\end{align}
By \eqref{22} and \eqref{23}, we get
\begin{align}
\sum_{n=0}^{\infty}a_{0,n}(\lambda)\frac{t^{n}}{n!}&=\overline{S_{\lambda}}(t)=e_{\lambda}^{1-\lambda}(t)S_{\lambda}(t)	\label{24}\\
&=\frac{d}{dt}\Big(e^{e_{\lambda}(t)-1}\Big)=\frac{d}{dt}\sum_{n=0}^{\infty}\phi_{n,\lambda}\frac{t^{n}}{n!}=\sum_{n=0}^{\infty}\phi_{n+1,\lambda}\frac{t^{n}}{n!}.\nonumber
\end{align}
Thus, by comparing the coefficients on both sides of \eqref{24}, we have
\begin{equation}
a_{0,n}(\lambda)=\phi_{n+1,\lambda},\quad (n\ge 0).\label{25}
\end{equation}
From \eqref{18} and \eqref{25}, we note that
\begin{align}
\phi_{n+1,\lambda}=a_{0,n}(\lambda)&=\sum_{k=0}^{n}\binom{n}{k}(1-\lambda)_{n-k,\lambda}a_{k,0}(\lambda)\label{26}\\
&=\sum_{k=0}^{n}\binom{n}{k}(1-\lambda)_{n-k,\lambda}\phi_{k,\lambda}. \nonumber	
\end{align}
Therefore, by \eqref{20}, \eqref{25} and \eqref{26}, we obtain the following theorem.
\begin{theorem}
For $n\ge 0$, we have
\begin{equation*}
\phi_{n+1,\lambda}=\sum_{k=0}^{n}\binom{n}{k}(1-\lambda)_{n-k,\lambda}\phi_{k,\lambda},
\end{equation*}
and
\begin{equation*}
\phi_{n,\lambda}=\sum_{k=0}^{n}\binom{n}{k}(-1)^{n-k}\langle 1-\lambda\rangle_{n-k,\lambda}\phi_{k+1,\lambda}.
\end{equation*}
\end{theorem}
Consider the initial degenerate sequence $\big(a_{n,0}(\lambda)\big)=\big(\phi_{n,\lambda}(x)\big)_{n\ge 0}$. Then the degenerate Euler-Seidel matrix corresponding to this sequence is given by
\begin{displaymath}
\left(\begin{matrix}
1 & x & (1-\lambda)x+x^{2} &  \cdots \\
1-\lambda+x & (2-\lambda)x+x^{2} & \cdots & \cdots \\
(1-\lambda)_{2,\lambda}+3(1-\lambda)x+x^{2} & \cdots & \cdots & \cdots \\
\vdots & \vdots & \vdots & \vdots
\end{matrix}\right).
\end{displaymath}
Note that
\begin{displaymath}
S_{\lambda}(t)=\sum_{n=0}^{\infty}a_{n,0}(\lambda)\frac{t^{n}}{n!}=\sum_{n=0}^{\infty}\phi_{n,\lambda}(x)\frac{t^{n}}{n!}=e^{x\big(e_{\lambda}(t)-1\big)}.
\end{displaymath}
From \eqref{22}, we have
\begin{align}
\sum_{n=0}^{\infty}a_{0,n}(\lambda)\frac{t^{n}}{n!}&=\overline{S}_{\lambda}(t)=e_{\lambda}^{1-\lambda}(t)S_{\lambda}(t)=e_{\lambda}^{1-\lambda}(t)e^{x(e_{\lambda}(t)-1)} \label{27}\\
&=\frac{1}{x}\frac{d}{dt}e^{x(e_{\lambda}(t)-1)}=\frac{1}{x}\sum_{n=0}^{\infty}\phi_{n+1,\lambda}(x)\frac{t^{n}}{n!}.\nonumber	
\end{align}
By \eqref{27}, we get
\begin{equation}
xa_{0,n}(\lambda)=\phi_{n+1,\lambda}(x),\quad (n\ge 0).\label{28}	
\end{equation}
Now, we observe that
\begin{align}
\sum_{n=0}^{\infty}a_{0,n}(\lambda)\frac{t^{n}}{n!}&=\overline{S_{\lambda}}(t)=e_{\lambda}^{1-\lambda}(t)S_{\lambda}(t) \label{29}\\
&=\sum_{n=0}^{\infty}\sum_{k=0}^{n}\binom{n}{k}(1-\lambda)_{n-k,\lambda}\phi_{k,\lambda}(x)\frac{t^{n}}{n!}. \nonumber	
\end{align}
Therefore, by \eqref{20}, \eqref{28} and \eqref{29}, we obtain the following theorem.
\begin{theorem}
For $n\ge 0$, we have
\begin{equation*}
\phi_{n+1,\lambda}(x)=x\sum_{k=0}^{n}\binom{n}{k}(1-\lambda)_{n-k,\lambda}\phi_{k,\lambda}(x),
\end{equation*}
and
\begin{displaymath}
x\phi_{n,\lambda}(x)=\sum_{k=0}^{n}\binom{n}{k}(\lambda-1)_{n-k,\lambda}\phi_{k+1,\lambda}(x).
\end{displaymath}
\end{theorem}
From \eqref{7}, we note that
\begin{align}
\sum_{n=0}^{\infty}\phi_{n,\lambda}^{\prime}(x)\frac{t^{n}}{n!}&=\frac{d}{dx}e^{x(e_{\lambda}(t)-1)}=\big(e_{\lambda}(t)-1\big)e^{x(e_{\lambda}(t)-1)}	\label{30}\\
&=e_{\lambda}(t) e^{x(e_{\lambda}(t)-1)}-e^{x(e_{\lambda}(t)-1)}\nonumber\\
&=e_{\lambda}^{\lambda}(t)e_{\lambda}^{1-\lambda}(t) e^{x(e_{\lambda}(t)-1)}-e^{x(e_{\lambda}(t)-1)} \nonumber\\
&=e_{\lambda}^{\lambda}(t)\overline{S_{\lambda}}(t)-S_{\lambda}(t). \nonumber
\end{align}
From \eqref{28} and \eqref{30}, we note that
\begin{align}
\sum_{n=0}^{\infty}\Big(\phi_{n,\lambda}^{\prime}(x)+\phi_{n,\lambda}(x)\Big)\frac{t^{n}}{n!}&=e_{\lambda}^{\lambda}(t)\overline{S_{\lambda}}(t)=e_{\lambda}^{\lambda}(t)\sum_{k=0}^{\infty}a_{0,k}(\lambda)\frac{t^{k}}{k!}\label{31}\\
&=\sum_{l=0}^{\infty}(\lambda)_{l,\lambda}\frac{t^{l}}{l!}\sum_{l=0}^{\infty}\frac{1}{x}\phi_{k+1,\lambda}(x)\frac{t^{k}}{k!}\nonumber\\
&=\sum_{n=0}^{\infty}\frac{1}{x}\sum_{k=0}^{n}\binom{n}{k}(\lambda)_{n-k,\lambda}\phi_{k+1,\lambda}(x)\frac{t^{n}}{n!}.\nonumber
\end{align}
Therefore, by \eqref{31}, we obtain the following theorem.
\begin{theorem}
For $n\ge 0$, we have
\begin{displaymath}
\sum_{k=0}^{n}\binom{n}{k}(\lambda)_{n-k,\lambda}\phi_{k+1,\lambda}(x)=x\Big(\phi_{n,\lambda}^{\prime}(x)+\phi_{n,\lambda}(x)\Big),
\end{displaymath}
where $\phi_{n,\lambda}^{\prime}(x)=\frac{d}{dx}\phi_{n,\lambda}(x)$.
\end{theorem}
Using \eqref{28} and \eqref{31}, we see that
\begin{align}
\frac{1}{x}\sum_{n=0}^{\infty}\phi_{n+1,\lambda}(x)\frac{t^{n}}{n!} &=\sum_{n=0}^{\infty}a_{0,n}(\lambda)\frac{t^{n}}{n!}=\overline{S_{\lambda}}(t)\label{32}\\
&=e_{\lambda}^{-\lambda}(t)\sum_{k=0}^{\infty}\Big(\phi_{k,\lambda}^{\prime}(x)+\phi_{k,\lambda}(x)\Big)\frac{t^{k}}{k!}\nonumber\\
&=\sum_{l=0}^{\infty}(-\lambda)_{l,\lambda}\frac{t^{l}}{l!}\sum_{k=0}^{\infty}\Big(\phi_{k,\lambda}^{\prime}(x)+\phi_{k,\lambda}(x)\Big)\frac{t^{k}}{k!}\nonumber\\
&=\sum_{n=0}^{\infty}\sum_{k=0}^{n}\binom{n}{k}(-\lambda)_{n-k,\lambda}\Big(\phi_{k,\lambda}^{\prime}(x)+\phi_{k,\lambda}(x)\Big)\frac{t^{n}}{n!}. \nonumber
\end{align}

Therefore, by \eqref{32}, we obtain the following theorem.
\begin{theorem}
For $n\ge 0$, we have
\begin{displaymath}
\phi_{n+1,\lambda}(x)=x\sum_{k=0}^{n}\binom{n}{k}(-1)^{n-k}\langle \lambda\rangle_{n-k,\lambda}\Big(\phi_{k,\lambda}^{\prime}(x)+\phi_{k,\lambda}(x)\Big).
\end{displaymath}
\end{theorem}
From Theorems 2.6 and 2.8, we note that
\begin{align}
x\phi_{n,\lambda}(x)&=\sum_{k=0}^{n}\binom{n}{k}(\lambda-1)_{n-k,\lambda}\phi_{k+1,\lambda}(x) \label{33}\\
&=x\sum_{k=0}^{n}\binom{n}{k}(\lambda-1)_{n-k,\lambda}\sum_{l=0}^{k}\binom{k}{l}(-\lambda)_{k-l,\lambda}\phi^{\prime}_{l,\lambda}(x)\nonumber\\
&\quad +x\sum_{k=0}^{n}\binom{n}{k}(\lambda-1)_{n-k,\lambda}\sum_{l=0}^{k}\binom{k}{l}(-\lambda)_{k-l,\lambda}\phi_{l,\lambda}(x).\nonumber
\end{align}
Now, using Lemma 1.1, we observe that
\begin{align}
&\sum_{k=0}^{n}\binom{n}{k}(\lambda-1)_{n-k,\lambda}\sum_{l=0}^{k}\binom{k}{l}(-\lambda)_{k-l,\lambda}\phi_{l,\lambda}^{\prime}(x) \label{34}\\
&=\sum_{l=0}^{n}\phi_{l,\lambda}^{\prime}(x)\sum_{k=l}^{n}\binom{n}{k}(\lambda-1)_{n-k,\lambda}\binom{k}{l}(-\lambda)_{k-l,\lambda} \nonumber \\
&=\sum_{l=0}^{n}\phi_{l,\lambda}^{\prime}(x)\binom{n}{l}(-1)_{n-l,\lambda}=\sum_{l=0}^{n}\phi_{l,\lambda}^{\prime}(x)\binom{n}{l}\langle 1\rangle_{n-l,\lambda}(-1)^{n-l}.\nonumber
\end{align}
By \eqref{33} and \eqref{34}, we get
\begin{align}
\phi_{n,\lambda}(x)&=\sum_{l=0}^{n}\binom{n}{l}\langle 1\rangle_{n-l,\lambda}(-1)^{n-l}\phi_{l,\lambda}^{\prime}(x)+\sum_{l=0}^{n}\binom{n}{l}\langle 1\rangle_{n-l,\lambda}(-1)^{n-l}\phi_{l,\lambda}(x)\label{35}\\
&=\sum_{l=0}^{n}\binom{n}{l}\langle 1\rangle_{n-l,\lambda}(-1)^{n-l}\phi_{l,\lambda}^{\prime}(x)+\sum_{l=0}^{n-1}\binom{n}{l}\langle 1\rangle_{n-l,\lambda}(-1)^{n-l}\phi_{l,\lambda}(x)+\phi_{n,\lambda}(x). \nonumber
\end{align}
Therefore, by \eqref{35}, we obtain the following theorem.
\begin{theorem}
For $n\ge 1$, we have
\begin{displaymath}
\sum_{k=0}^{n-1}\binom{n}{k}\langle 1\rangle_{n-k,\lambda}(-1)^{k}\phi_{k,\lambda}(x)=\sum_{k=0}^{n}\binom{n}{k}\langle 1\rangle_{n-k,\lambda}(-1)^{k-1}\phi_{k,\lambda}^{\prime}(x).
\end{displaymath}
\end{theorem}
Consider the initial degenerate sequence $(a_{n,0}(\lambda))_{n\ge 0}=(F_{n,\lambda})_{n\ge 0}$, (see \eqref{10}). Then degenerate Euler-Seidel matrix corresponding to this sequence is given by
\begin{displaymath}
\left(\begin{matrix}
1 & 1 & 3-\lambda &\cdots  \\
2-\lambda & 4-\lambda & 16-7 \lambda+11 \lambda^2 & \cdots  \\
6-6\lambda+2\lambda^{2} & 20-12\lambda+12\lambda^{2} & 104-32 \lambda+38 \lambda^2 +18 \lambda^3 &\cdots  \\
\vdots & \vdots &\vdots & \vdots \\
\end{matrix}\right).
\end{displaymath}
Note that
\begin{equation}
S_{\lambda}(t)=\sum_{n=0}^{\infty}a_{n,0}(\lambda)\frac{t^{n}}{n!}=\sum_{n=0}^{\infty}F_{n,\lambda}\frac{t^{n}}{n!}=\frac{1}{2-e_{\lambda}(t)}. \label{36}
\end{equation}
From \eqref{22} and \eqref{36}, we have
\begin{align}
\sum_{n=0}^{\infty}a_{0,n}(\lambda)\frac{t^{n}}{n!}&=\overline{S_{\lambda}}(t)=e_{\lambda}^{1-\lambda}(t)S_{\lambda}(t)=e_{\lambda}^{-\lambda}(t)\frac{e_{\lambda}(t)-2+2}{2-e_{\lambda}(t)} \label{37} \\
&=-e_{\lambda}^{-\lambda}(t)+\frac{2}{2-e_{\lambda}(t)}e_{\lambda}^{-\lambda}(t)\nonumber\\
&=\sum_{n=0}^{\infty}\bigg(2\sum_{k=0}^{n}\binom{n}{k}(-\lambda)_{n-k,\lambda}F_{k,\lambda}-(-\lambda)_{n,\lambda}\bigg)\frac{t^{n}}{n!}.\nonumber
\end{align}
For $n\ge 1$, by \eqref{18} and \eqref{37}, we get
\begin{equation}
2\sum_{k=0}^{n}\binom{n}{k}(-\lambda)_{n-k,\lambda}F_{k,\lambda}-(-\lambda)_{n,\lambda}=a_{0,n}(\lambda)=\sum_{k=0}^{n}\binom{n}{k}(1-\lambda)_{n-k,\lambda}F_{k,\lambda}. \label{38}	
\end{equation}
From \eqref{20}, \eqref{38} and Lemma 1.1, we note that
\begin{align}
F_{n,\lambda}&=a_{n,0}(\lambda)=\sum_{k=0}^{n}\binom{n}{k}(-1)^{n-k}\langle 1-\lambda\rangle_{n-k,\lambda}a_{0,k}(\lambda) \label{39} \\
&=\sum_{k=0}^{n}\binom{n}{k}(\lambda-1)_{n-k,\lambda}\bigg(2\sum_{l=0}^{k}\binom{k}{l}(-\lambda)_{k-l,\lambda}F_{l,\lambda}-(-\lambda)_{k,\lambda}\bigg) \nonumber\\
&=2\sum_{l=0}^{n}F_{l,\lambda}\sum_{k=l}^{n}\binom{n}{k}(\lambda-1)_{n-k,\lambda}\binom{k}{l}(-\lambda)_{k-l,\lambda}\nonumber\\
&\quad -\sum_{k=0}^{n}\binom{n}{k}(\lambda-1)_{n-k,\lambda}(-\lambda)_{k,\lambda}\nonumber\\
&=2\sum_{l=0}^{n}\binom{n}{l}F_{l,\lambda}(-1)_{n-l,\lambda}-(-1)_{n,\lambda}\nonumber\\
&=2\sum_{l=1}^{n}\binom{n}{l}F_{l,\lambda}(-1)^{n-l}\langle 1\rangle_{n-l,\lambda}+(-1)^{n-1}\langle 1\rangle_{n,\lambda}. \nonumber
\end{align}
Therefore, by \eqref{38} and \eqref{39}, we obtain the following theorem.
\begin{theorem}
For $n\ge 1$, we have
\begin{displaymath}
F_{n,\lambda}= 2\sum_{l=1}^{n}\binom{n}{l}F_{l,\lambda}(-1)^{n-l}\langle 1\rangle_{n-l,\lambda}+(-1)^{n-1}\langle 1\rangle_{n,\lambda},
\end{displaymath}
and
\begin{displaymath}
\sum_{k=0}^{n}\binom{n}{k}(1-\lambda)_{n-k,\lambda}F_{k,\lambda}=2\sum_{k=0}^{n}\binom{n}{k}(-1)^{n-k}\langle \lambda\rangle_{n-k,\lambda}F_{k,\lambda}+(-1)^{n-1}\langle \lambda\rangle_{n,\lambda}.
\end{displaymath}
\end{theorem}
Consider the initial degenerate sequence $\big(a_{n,0}(\lambda)\big)_{n\ge 0}=\big(F_{n,\lambda}(x)\big)_{n\ge 0}$. Then, by \eqref{10}, we get
\begin{align}
S_{\lambda}(t)=\sum_{n=0}^{\infty}a_{n,0}(\lambda)\frac{t^{n}}{n!}=\sum_{n=0}^{\infty}F_{n,\lambda}(x)\frac{t^{n}}{n!}=\frac{1}{1-x(e_{\lambda}(t)-1)}.\label{40}
\end{align}
From \eqref{22} and \eqref{40}, we have
\begin{equation}
\sum_{n=0}^{\infty}a_{0,n}(\lambda)\frac{t^{n}}{n!}=\overline{S}_{\lambda}(t)=e_{\lambda}^{1-\lambda}(t)S_{\lambda}(t)=e_{\lambda}^{1-\lambda}(t)\frac{1}{1-x(e_{\lambda}(t)-1)}. \label{41}
\end{equation}
Now, we observe that
\begin{align}
\frac{d}{dt}S_{\lambda}(t)&=\frac{d}{dt}\bigg(\frac{1}{1-x(e_{\lambda}(t)-1)}\bigg)=\frac{xe_{\lambda}^{1-\lambda}(t)}{\big(1-x(e_{\lambda}(t)-1)\big)^{2}}.\label{42}
\end{align}
Thus, by \eqref{41} and \eqref{42}, we get
\begin{align}
\sum_{n=0}^{\infty}a_{0,n}(\lambda)\frac{t^{n}}{n!}&=\overline{S}_{\lambda}(t)=\bigg(\frac{1}{x}-\big(e_{\lambda}(t)-1\big)\bigg)\frac{d}{dt}S_{\lambda}(t)\label{43} \\
&=\bigg(\frac{1}{x}-e_{\lambda}(t)+1\bigg)\sum_{k=0}^{\infty}F_{k+1,\lambda}(x)\frac{t^{k}}{k!}\nonumber\\
&=\sum_{n=0}^{\infty}\bigg(\frac{1}{x}F_{n+1,\lambda}(x)-\sum_{k=0}^{n}\binom{n}{k}F_{k+1,\lambda}(x)(1)_{n-k,\lambda}+F_{n+1,\lambda}(x)\bigg)\frac{t^{n}}{n!}. \nonumber
\end{align}
Comparing the coefficients on both sides of \eqref{43}, we have
\begin{align}
a_{0,n}(\lambda)&=\frac{1}{x}F_{n+1,\lambda}(x)-\sum_{k=0}^{n}\binom{n}{k}F_{k+1,\lambda}(x)(1)_{n-k,\lambda}+F_{n+1,\lambda}(x) \label{44}\\
&=\frac{1}{x}F_{n+1,\lambda}(x)-\sum_{k=0}^{n-1}\binom{n}{k}F_{k+1,\lambda}(x)(1)_{n-k,\lambda}\nonumber\\
&=\frac{1}{x}F_{n+1,\lambda}(x)-\sum_{k=1}^{n}\binom{n}{k-1}F_{k,\lambda}(x)(1)_{n-k+1,\lambda}\nonumber\\
&=\frac{1}{x}F_{n+1,\lambda}(x)-\sum_{k=1}^{n}\binom{n}{k-1}F_{k,\lambda}(x)(1-\lambda)_{n-k,\lambda}.\nonumber
\end{align}
By \eqref{18}, we get
\begin{equation}
a_{0,n}(\lambda)=\sum_{k=0}^{n}\binom{n}{k}(1-\lambda)_{n-k,\lambda}F_{k,\lambda}(x).\label{45}
\end{equation}
From \eqref{44} and \eqref{45}, we note that
\begin{equation}
F_{n+1,\lambda}(x)-x\sum_{k=1}^{n}\binom{n}{k-1}F_{k,\lambda}(x)(1-\lambda)_{n-k,\lambda}=x\sum_{k=0}^{n}\binom{n}{k}(1-\lambda)_{n-k,\lambda}F_{k,\lambda}(x).\label{46}
\end{equation}
Thus, by \eqref{46}, we get
\begin{align}
F_{n+1,\lambda}(x)&=x\sum_{k=1}^{n}\bigg(\binom{n}{k-1}+\binom{n}{k}\bigg)(1-\lambda)_{n-k,\lambda}F_{k,\lambda}(x)+(1-\lambda)_{n,\lambda}x \label{47}\\
&=x\sum_{k=1}^{n}\binom{n+1}{k}(1-\lambda)_{n-k,\lambda}F_{k,\lambda}(x)+(1-\lambda)_{n,\lambda}x\nonumber\\
&=x\sum_{k=0}^{n}\binom{n+1}{k}(1-\lambda)_{n-k,\lambda}F_{k,\lambda}(x). \nonumber	
\end{align}
Therefore, by \eqref{47}, we obtain the following theorem.
\begin{theorem}
For $n\ge 1$, we have
\begin{displaymath}
F_{n,\lambda}(x)=x\sum_{k=0}^{n-1}\binom{n}{k}(1-\lambda)_{n-1-k,\lambda}F_{k,\lambda}(x).
\end{displaymath}
\end{theorem}
From \eqref{47}, we note that
\begin{align}
F_{n+1,\lambda}(x)&=x\sum_{k=0}^{n}\binom{n+1}{k}(1-\lambda)_{n-k,\lambda}F_{k,\lambda}(x)\label{48}\\
&=x\sum_{k=0}^{n}\bigg(\binom{n}{k}+\binom{n}{k-1}\bigg)(1-\lambda)_{n-k,\lambda}F_{k,\lambda}(x) \nonumber\\
&=x\sum_{k=0}^{n}\binom{n}{k}(1-\lambda)_{n-k,\lambda}F_{k,\lambda}(x)+x\sum_{k=1}^{n}\binom{n}{k-1}(1-\lambda)_{n-k,\lambda}F_{k,\lambda}(x)\nonumber\\
&=x\sum_{k=0}^{n}\binom{n}{k}(1-\lambda)_{n-k,\lambda}F_{k,\lambda}(x)+x\sum_{k=0}^{n-1}\binom{n}{k}(1-\lambda)_{n-k-1,\lambda}F_{k+1,\lambda}(x)\nonumber\\
&=x\sum_{k=0}^{n}\binom{n}{k}(1-\lambda)_{n-k,\lambda}F_{k,\lambda}(x)+x\sum_{k=0}^{n}\binom{n}{k}(1)_{n-k,\lambda}F_{k+1,\lambda}(x)-xF_{n+1,\lambda}(x). \nonumber
\end{align}
Therefore, by \eqref{48}, we obtain the following theorem.
\begin{theorem}
For $n\ge 0$, we have
\begin{displaymath}
F_{n+1,\lambda}(x)=\frac{x}{1+x}\sum_{k=0}^{n}\binom{n}{k}\Big((1-\lambda)_{n-k,\lambda}F_{k,\lambda}(x)+(1)_{n-k,\lambda}F_{k+1,\lambda}(x)\Big).
\end{displaymath}
\end{theorem}
We note from \eqref{6} and \eqref{9} that
\begin{equation}
F_{n+1,\lambda}(x)=\sum_{k=0}^{n+1} {n+1 \brace k}_{\lambda} x^{k} \int_{0}^{\infty}y^{k}e^{-y} dy=\int_{0}^{\infty} \phi_{n+1,\lambda}(xy)e^{-y} dy. \label{49}
\end{equation}
Thus, by \eqref{6}, \eqref{49} and Theorem 2.6, we get
\begin{align}
F_{n+1,\lambda}(x)&=\int_{0}^{\infty}\phi_{n+1,\lambda}(xy)e^{-y}dy=x\sum_{k=0}^{n}	\binom{n}{k}(1-\lambda)_{n-k,\lambda}\int_{0}^{1}y\phi_{k,\lambda}(xy)e^{-y}dy \label{50} \\
&=x\sum_{k=0}^{n}\binom{n}{k}(1-\lambda)_{n-k,\lambda}\sum_{j=0}^{k}{k \brace j}_{\lambda}x^{j}\int_{0}^{\infty}y^{j+1}e^{-y}dy \nonumber\\
&=x\sum_{k=0}^{n}\binom{n}{k}(1-\lambda)_{n-k,\lambda}\sum_{j=0}^{k}{k \brace j}_{\lambda}x^{j}(j+1)!.\nonumber
\end{align}
Note that
\begin{align}
\frac{d}{dx}\Big(xF_{k,\lambda}(x)\Big)&=\frac{d}{dx}\bigg(\sum_{j=0}^{k}{k \brace j}_{\lambda}j!x^{j+1}\bigg)\label{51} \\
&=\sum_{j=0}^{k}{k \brace j}_{\lambda}(j+1)!x^{j}.\nonumber
\end{align}
From \eqref{50} and \eqref{51}, we have
\begin{align}
F_{n+1,\lambda}(x)&=x\sum_{k=0}^{n}\binom{n}{k}(1-\lambda)_{n-k,\lambda}\sum_{j=0}^{k}{k \brace j}_{\lambda}x^{j}(j+1)!\label{52}\\
&=x\sum_{k=0}^{n}\binom{n}{k}(1-\lambda)_{n-k,\lambda}\frac{d}{dx}\Big(xF_{k,\lambda}(x)\Big)\nonumber\\
&=x\sum_{k=0}^{n}\binom{n}{k}(1-\lambda)_{n-k,\lambda}\Big(F_{k,\lambda}(x)+xF_{k,\lambda}^{\prime}(x)\Big). \nonumber
\end{align}
Therefore, by \eqref{52}, we obtain the following theorem.
\begin{theorem}
For $n\ge 0$, we have
\begin{equation}
F_{n+1,\lambda}(x)=x\sum_{k=0}^{n}\binom{n}{k}(1-\lambda)_{n-k,\lambda}\Big(F_{k,\lambda}(x)+xF_{k,\lambda}^{\prime}(x)\Big).\label{53}
\end{equation}
\end{theorem}
From Theorem 2.11 and \eqref{53}, we note that
\begin{align}
&x\sum_{k=0}^{n}\binom{n+1}{k}(1-\lambda)_{n-k,\lambda}F_{k,\lambda}(x)=F_{n+1,\lambda}(x)\label{54}\\
&=x\sum_{k=0}^{n}\binom{n}{k}(1-\lambda)_{n-k,\lambda}F_{k,\lambda}(x)+x\sum_{k=0}^{n}\binom{n}{k}(1-\lambda)_{n-k,\lambda}xF_{k,\lambda}^{\prime}(x).\nonumber
\end{align}
Thus, by \eqref{54}, we get
\begin{align}
&\sum_{k=0}^{n}\binom{n}{k}(1-\lambda)_{n-k,\lambda}xF_{k,\lambda}^{\prime}(x) \label{55}\\
&=\sum_{k=0}^{n}\bigg(\binom{n+1}{k}-\binom{n}{k}\bigg)(1-\lambda)_{n-k,\lambda}F_{k,\lambda}(x)\nonumber\\
&=\sum_{k=1}^{n}\binom{n}{k-1}(1-\lambda)_{n-k,\lambda}F_{k,\lambda}(x). \nonumber	
\end{align}
Therefore, by \eqref{55}, we obtain the following theorem.
\begin{theorem}
For $n\ge 1$, we have
\begin{displaymath}
\sum_{k=1}^{n}\binom{n}{k-1}(1-\lambda)_{n-k,\lambda}F_{k,\lambda}(x)= \sum_{k=0}^{n}\binom{n}{k}(1-\lambda)_{n-k,\lambda}xF_{k,\lambda}^{\prime}(x).
\end{displaymath}
\end{theorem}
\section{Conclusion}

In this paper, we successfully developed and analyzed a degenerate Euler-Seidel matrix method by introducing a parameter $\lambda$ into the fundamental recurrence relation. This generalized approach not only preserves the structure of the classical Euler-Seidel method as the parameter $\lambda \to 0$ but also provides a powerful framework for studying degenerate versions of combinatorial sequences. We established the key transformation formulas relating the initial degenerate sequence $\big(a_{n,0}(\lambda)\big)$ to the final degenerate sequence $\big(a_{0,n}(\lambda)\big)$, demonstrating a direct generalization of the classical binomial identities using the $\lambda$-generalized falling and rising factorials, $(1-\lambda)_{n-k,\lambda}$ and $\langle1-\lambda \rangle_{n-k,\lambda}$. Most significantly, we derived the degenerate Seidel's formula for the exponential generating functions:$$\overline{S_{\lambda}}(t)=e_{\lambda}^{1-\lambda}(t)S_{\lambda}(t).$$The applications of this degenerate method to sequences such as the degenerate Bell numbers and polynomials, and the degenerate Fubini numbers and polynomials, yielded a rich set of novel combinatorial identities, showing the utility and broad applicability of our generalization. This work opens avenues for future research into other generalized and degenerate matrix methods and their connections to special numbers and polynomials in combinatorial analysis.

\end{document}